\documentclass{article}
\usepackage{graphicx}
\usepackage{amsmath}
\usepackage{epsf}
\usepackage{graphpap}
\usepackage{moreverb}
\usepackage{placeins}
\usepackage{float}
\usepackage{placeins}
\usepackage{setspace}
\usepackage{caption}
\usepackage{csquotes}

\title{Cost-Effectiveness Analysis of COVID-19 Vaccination : A
review of some Vaccination Models}
\author{Rehana Naz  $^{a,*}$, Andrew Omame $^b$, Mariano
Torrisi$^{c}$  \\
\centerline{$^{a}$Department of Mathematics,}\\
   \centerline{ Lahore School of Economics}\\
   \centerline{ Lahore, 53200, Pakistan}\\
   \centerline{$^b$ Department of Mathematics  }\\
   \centerline{  Federal University of Technology}\\
    \centerline{Owerri, Nigeria  }\\
     \centerline{$^c$ Dipartimento di Matematica e Informatica  }\\
   \centerline{  UniversitA degli Studi di Catania}\\
    \centerline{Viale Andrea Doria, 6, Catania, 95125, Italy.   }\\
    $^*$ Corresponding author emails: rehananaz\_qau@yahoo.com, \\drrehana@lahoreschool.edu.pk }

\begin{document}

\maketitle

\begin{abstract}
The sudden and rapid spread of the COVID\_19 pandemic with its
terrible consequences has put the management of governments and the
various world institutions into a crisis. They have been subjected
to a considerable economic effort to be taken to combat the spread
of the pandemic. The economic investment for the research and
purchase of vaccines intended for populations is subject to
cos-benefit analyzes in various situations in different cases. In
this review work, several recent models are analyzed where the
appearance of the components is coupled with the economic aspect.
The analysis of these models is detailed and the results discussed
from different points of view.
\end{abstract}
\section*{Keywords}Cost-Effectiveness Analysis, COVID-19 Vaccination models, Environmental transmission, Co-circulation of COVID-19 and Arboviruses
\section{ Introduction}

The emergence of pandemics has revealed the vulnerability of mankind
to this phenomenon. COVID-19 has rapidly spread across the globe,
resulting in both fatalities and substantial social problems. The
necessity to intervene to curb and control its progression has
initiated significant scientific and political debates. Government
initiatives not only impose restrictions on the daily life of the
population but also result in substantial financial costs. The most
widespread measure adopted almost everywhere was vaccination.
Obtaining available vaccines from the market (and considering their
cost) posed an issue. In fact, procuring sufficient quantities of
vaccines to adequately cater to the population was not an easy task.
Limited accessibility has posed numerous challenges for the
governments of the least developed nations. We often resorted to
cost-benefit analyses to understand the correct allocation of
resources. Since, in general, we are dealing with a priori
simulations, the development of mathematical models on pandemics has
been of great help. Particularly focusing on models that incorporate
economic variables into the system of equations, there is a limited
amount of published material on these combined themes. Our aim in
this work is to provide the reader with a brief but comprehensive
overview. We have concentrated on several mathematical models that
not only track the evolution of variables describing pandemics under
the influence of vaccination but also incorporate a
cost-effectiveness analysis.     Naturally, our coverage may not be
exhaustive given the dynamic nature of research boundaries.

  \par
Gold \cite{ce} is a foundational reference book providing guidelines
and methodologies specific to cost-effectiveness analysis in the
health and medicine domain. Muennig et al. \cite{econ} offer a
practical guide for those involved in health economics, providing
insights into the application of cost-effectiveness analysis in
health-related decision-making. Levin et al. \cite{edu} provide a
broader perspective on the application of cost-effectiveness
analysis beyond healthcare; this book explores its use in education
and social programs. Two distinct methods, specifically the average
cost-effectiveness ratio (ACER) and the incremental
cost-effectiveness ratio (ICER), play a pivotal role in determining
the cost-effective intervention plan among the various combinations
under consideration. ACER  focuses solely on comparing the single
intervention technique to its baseline alternative. This involves
assessing the proportion of the overall cost of the intervention to
the total number of infections it prevents.  The ACER formula, as
presented by \cite{Ref25,Ref26}, is expressed as:

    \begin{equation}\label{ACER}
    \text {ACER }=\frac{\text {Total cost produced by intervention }}{\text {Total number of infection averted}}.
    \end{equation}
    The denominator of equation \eqref{ACER} is derived by subtracting the number of infected people under control from the number of infectious people not under control.

    On the other hand, the ICER is concerned with comparing the differences in costs and health outcomes between two alternative intervention strategies that compete for the same resources.  It is the ratio of the increase or decrease in the total number of infections prevented by two different techniques to the change in the expenses of those two treatments. Simply, ICER is determined as follows:

\begin{equation}
\text { ICER }=\frac{\text { Change in intervention costs }}{\text {
Change in total number of infection averted }}.
\end{equation}

We have focused on several mathematical models that not only track
the evolution of variables describing pandemics under the influence
of vaccination but also incorporate a cost-effectiveness analysis.
Naturally, our coverage of mathematical models may not be
exhaustive, given the dynamic nature of research boundaries.In their
respective studies, Li et al. \cite{ce1} concentrated on diabetes
mellitus, Otieno et al. \cite{ce2} addressed malaria control
strategies, Myran et al. \cite{ce3} explored hepatitis B virus, and
Fesenfeld \cite{ce4} investigated human papillomavirus vaccination,
each performing a comprehensive cost-effectiveness analysis in their
investigations. Jit and Mibei \cite{ce5} conducted a systematic
review that discusses the cost-effectiveness of vaccination
programs, providing valuable insights into the economic value of
preventive health measures. In \cite{ce6, ce6b}  an in-depth
exploration of the theoretical foundations of cost-effectiveness
analysis, established a conceptual framework for understanding its
underlying principles. Brent \cite{ce7} contributed to the field by
exploring the Societal Perspective in Healthcare. The work of
\cite{ce8} specifically examined COVID-19 Inactivated Vaccines.
References \cite{Ref5}-\cite{Ref21} investigate various case
studies, both empirical and theoretical, discussing COVID-19 from
cost-effectiveness perspectives, and we provide a detailed
discussion of these studies in our article.

 Here we tried to offer essential information on as many works as possible within the review field. In our opinion, comprehensive and
  exhaustive information on the contents of the works holds greater significance than the opinions of
  the review authors provided through in-depth analysis. The plan of the paper is outlined as follows: In Section
\ref{Cost-effectiveness}, a detailed review of the
cost-effectiveness of COVID-19 vaccination models is presented. The
discussions and future directions are provided in Section 3.

\section{Cost-effectiveness of COVID-19 vaccination models}\label{Cost-effectiveness}
We conducted an extensive literature search by exploring the
prominent scholarly electronic database PubMed/MEDLINE/Scopus. The
keywords used in the search included: \enquote{COVID-19},
\enquote{SIR}, \enquote{Epidemiological models},
\enquote{Cost-effectiveness}, and \enquote{vaccine}. Additionally,
we employed \enquote{Medical Subject Headings} (MeSH) terms and
truncated words when necessary. The search was confined to studies
specifically conducted on COVID-19, excluding investigations on
other diseases.

\subsection{SEAIR framework with environmental transmission of COVID-19}

%Optimal control and comprehensive cost-effectiveness analysis for COVID-19}
Using data from the Kingdom of Saudi Arabia as a case study, Asamoah
\textit{et al.} \cite{Ref5} investigated a COVID-19 model with
optimal control and thorough cost-effectiveness analysis. The model
took into account four time-dependent control measures, including
social segregation protocols, maintaining personal hygiene by using
alcohol-based detergents to clean contaminated surfaces,
implementing appropriate and secure measures for individuals who are
exposed, asymptomatic, or symptomatic with infection,and carry out
fumigation in educational institutions across all levels, sports
facilities, commercial spaces, and places of worship. The system of
equations proposed by Asamoah \textit{et al.} \cite{Ref5}  is as
follows:

\begin{equation}\label{Ref5}
\begin{split}
\frac{d S}{d t} & =\Lambda-\left(\beta_1 E+\beta_2 I+\beta_3 A\right) \frac{S}{N}-\beta_4 B \frac{S}{N}-d S, \\
\frac{d E}{d t} & =\left(\beta_1 E+\beta_2 I+\beta_3 A\right) \frac{S}{N}+\beta_4 B \frac{S}{N}-(d+\delta) E, \\
\frac{d I}{d t} & =\delta(1-\tau)  E-\left(d+d_1+\gamma_1\right) I, \\
\frac{d A}{d t} & =\tau \delta E-\left(d+\gamma_2\right) A, \\
\frac{d R}{d t} & =\gamma_1 I+\gamma_2 A-d R, \\
\frac{d B}{d t} & =\psi_1 E+\psi_2 I+\psi_3 A-\phi B,
\end{split}
\end{equation}
where, the model variables $S$, $E$,  $A$, $I$, $R$ are Susceptible,
Exposed, Asymptomatic, Symptomatic and Recovered individuals
respectively with $N$  as the total population, $B(t)$ represents
the assumed concentration of the SARS-CoV-2 in the environment. The
constitutive parameters are defined in Table \ref{Ref5Table} in the
Appendix:

By conducting a thorough cost-effectiveness analysis for various
control scenarios, numerical simulations were performed. The results
revealed that implementing protocols for physical or social
distancing is the most effective and economically viable
intervention in the Kingdom of Saudi Arabia, following a
comprehensive comparison of various control measures.

In a separate study, Asamoah \textit{et al.} \cite{Ref6} examined
the global stability and economic viability of a COVID-19 model
while taking environmental transmission into account. They employed
real data from Ghana in their study. The most  economically
efficient way  to curb the transmission of COVID-19 in Ghana,
according to their study, is to combine safety measures like
practicing appropriate coughing etiquette, encouraging the use of
disposable tissues or clothing to cover coughs and sneezes and
emphasize the importance of hand-washing after such instances. In
reference \cite{Ref6}, the model presented below was studied
\begin{equation}\label{Ref6}
\begin{split}
& \frac{d S}{d t}=\Lambda-\omega S-\beta(I S+\eta A S)-\beta_1 V S+\rho R, \\
& \frac{d E}{d t}=\beta(I S+\eta A S)+\beta_1 V S-k_2(1-\gamma) E-k_1 \gamma E-\omega E, \\
& \frac{d A}{d t}=k_2(1-\gamma) E-\omega A-v_1 \phi A-v_2(1-\phi) A, \\
& \frac{d I}{d t}=k_1 \gamma E+v_1 \phi A-\epsilon I-(\omega+\alpha) I, \\
& \frac{d R}{d t}=v_2(1-\phi) A+\epsilon I-\rho R-\omega R, \\
& \frac{d V}{d t}=m_1 A+m_2 I-\tau_1 V,
\end{split}
\end{equation}
where, the model variables $S$, $E$, $A$, $I$ and $R$ represent the
susceptible, exposed/pre asymptomatic, Asymptomatic,  Symptomatic
and Recovered individuals. The variable $V$ is Virus in the
Environment. The model parameters are defined in Table
\ref{Ref6Table} in the Appendix section:

\subsection{SEIAHRD Model for COVID-19 Analysis in Nigeria: Incorporating Real Data}

%Mathematical modelling and optimal cost-effective control of COVID-19 transmission dynamics}
Olaniyi et al \cite{Ref7} proposed a mathematical framework to
examine the dynamics of COVID-19 transmission, with a specific
emphasis on optimal and cost-effective control measures. They used
actual COVID-19 data from Nigeria to fit their model. They found
that in terms of reducing the burden of COVID-19, the best
preventive measure-which includes immunisation, the use of face
masks and hand sanitizers performed significantly better than
management control. Management control includes providing adequate
care for patients who are hospitalised in order to ensure a speedy
recovery and prevent deaths from complications.  They found that
implementing both control measures together is the most
cost-effective approach when compared to the individual
implementation of each control measure.   Olaniyi et al \cite{Ref7}
considered the following model:
    \begin{equation}\label{Ref7}
    \begin{split}
    \frac{d S}{d t} & =-\frac{\beta S\left(I+\varepsilon_1 A+\varepsilon_2 H\right)}{N-D} \\
    \frac{d E}{d t} & =\frac{\beta S\left(I+\varepsilon_1 A+\varepsilon_2 H\right)}{N-D}-\alpha E\\
    \frac{d I}{d t} & =l_1 \alpha E-\left(h_1+r_1+\delta_1\right) I \\
    \frac{d A}{d t} & =\left(1-l_1\right) \alpha E-\left(h_2+r_2\right) A \\
    \frac{d H}{d t} & =h_1 I+h_2 A-\left(\gamma+\delta_2\right) H \\
    \frac{d R}{d t} & =r_1 I+r_2 A+\gamma H \\
    \frac{d D}{d t} & =\delta_1 I+\delta_2 H \\
    \end{split}
    \end{equation}
    where, the model variables $S$, $E$, $A$, $I$, $R$, $H$ and $D$ represent the susceptible, exposed, Asymptomatic,  Symptomatic, Recovered, Hospitalized, Death
classes.  Definitions of constitutive parameters can be found in
Table \ref{Ref7Table} of the Appendix section.

\subsection{$SEAICC_WHR$ Model: Understanding COVID-19 Dynamics in Peru}
%Cost-effectiveness of a mathematical modeling with optimal control approach of spread of COVID-19 pandemic: A case study in Peru}

The system of differential equations governing the model proposed by
Kouidere et al.\cite{Ref8}  is as follows:

        \begin{equation}\label{Ref8}
    \left\{\begin{array}{l}
    \frac{d S(t)}{d t}=\Lambda-\mu S(t)-\beta_1 \frac{S(t) A(t)}{N}-\beta_2 \frac{S(t) I(t)}{N} \\
    \frac{d E(t)}{d t}=\beta_1 \frac{S(t) A(t)}{N}+\beta_2 \frac{S(t)](t)}{N}-\left(\mu+\alpha_1+\alpha_2\right) E(t) \\
    \frac{d A(t)}{d t}=\alpha_1 E(t)-\left(\theta_1+\chi+\mu\right) A(t) \\
    \frac{d I(t)}{d t}=\alpha_2 E(t)+\chi A(t)-\left(\gamma+\theta_2+\mu\right) I(t) \\
    \frac{d C(t)}{d t}=(1-\varrho) \gamma I(t)-\left(\eta_1+\mu+\delta_1\right) C(t) \\
    \frac{d C_W(t)}{d t}=\varrho \gamma I(t)-\left(\eta_2+\mu+\delta_2\right) C_W(t) \\
    \frac{d H(t)}{d t}=\theta_1 I(t)+\theta_2 I(t)+\eta_1 C(t)+\eta_2 C_W(t)-\left(\mu+\sigma+\delta_3\right) H(t) \\
    \frac{d R(t)}{d t}=\sigma H(t)-\mu R(t)
    \end{array}\right.
    \end{equation}

The model variables $S$, $E$, $A$, $I$, and $R$ represent
susceptible, exposed, asymptomatic, symptomatic infected, and
recovered individuals, respectively. The compartment $C$ denotes
those infected with complications, while $C_W$  represents
individuals with complications and chronic diseases. Individuals
under lockdown in hospitals with follow-up and health monitoring are
grouped as $H$. The model parameters are defined in Table
\ref{Ref8Table} in the Appendix Section.

   Kouidere \textit{et al.} \cite{Ref8} carried out a cost-effectiveness analysis of a COVID-19 mathematical model based on data from Peru. Utilizing Pontryagin's Maximum Principle, they demonstrated that the most effective strategy for controlling the spread of COVID-19 in Peru involves a combination of conducting awareness campaigns and implementing quarantine alongside treatment for infected individuals.
\subsection{A $SEIAQHRM$ model: Self-Protection, Vaccination, and Disinfectant Strategies for COVID-19 Control}
%A Mathematical Evaluation of the Cost-Effectiveness of Self-Protection, Vaccination, and Disinfectant Spraying for COVID-19 Control}

The system of nonlinear ordinary differential equations for the
model proposed by Nana-Kyere and Seidu \cite{Ref9} is as follows:
\begin{equation}\label{Ref9}
\begin{split}
& \frac{d S}{d t}=\Lambda-\mu S-\eta_2 M S-\eta_1 \frac{\text { IS }}{N}-\frac{\psi A S}{N} \eta_1, \\
& \frac{d E}{d t}=\eta_2 M S+\eta_1 \frac{\text { IS }}{N}+\frac{\psi A S}{N} \eta_1-\left((1-\theta) \omega+\theta \rho+\delta_1+\mu\right) E, \\
& \frac{d I}{\mathrm{~d} t}=(1-\theta) \omega E-\left(\tau_1+\xi_1+\gamma+\mu\right) I, \\
& \frac{d A}{d t}=\theta \rho E-\left(\tau_2+\mu\right) A, \\
& \frac{d Q}{d t}=\delta_1 E-\left(\phi_1+\delta_2+\mu\right) Q, \\
& \frac{d H}{d t}=\gamma I+\delta_2 Q-\left(\phi_2+\xi_2+\mu\right) H, \\
& \frac{d R}{d t}=\tau_1 I+\tau_2 A+\phi_1 Q+\phi_2 H-\mu R, \\
& \frac{d M}{d t}=q_1 I+q_2 A-q_3 M,
\end{split}
\end{equation}
where, the model variables are defined as follows: $S$: Susceptible,
$E$: Exposed, $I$: Infected, $A$: asymptomatically infected, $Q$:
Quarantined,  $H$: Hospitalized, $R$: Recovered, and $M$: virus
concentration in the environment. The model parameters are defined
in Table \ref{Ref9Table} in the Appendix section.

Nana-Kyere and Seidu \cite{Ref9} investigated a mathematical model
to understand the dynamics of COVID-19 transmission. They conducted
a cost-effectiveness evaluation of control measures, including
self-protection, vaccination, and disinfectant spraying. Their study
revealed that the most cost-effective approach to control the
transmission of the virus involves the combined implementation of
self-protection and environmental control measures.

\subsection{The  SEIHR framework  to Model COVID-19 transmission Dynamics in Italy and Spain}

Srivastav et al. \cite{Ref11} developed and analyzed a compartmental
COVID-19 model that categorizes the population into old and young
groups. By computing both the disease-free equilibrium and the basic
reproduction number, the study provides a comprehensive
understanding of the impact of COVID-19. The estimation of key
parameters using real-life data from Italy and Spain, employing the
least square method, enhances the empirical validity of the model.
By emphasizing parameters with a significant effect on the basic
reproduction number, the sensitivity analysis provides insight into
essential factors shaping transmission dynamics. To investigate the
most efficient, time-dependent, and cost-effective control
strategies that can lower the number of infections within a given
time frame, the model is expanded to include optimal control
problems. They characterized the controls as follows: mitigating the
person-to-person transmission of disease through measures like
social distancing, awareness campaigns, and enhanced sanitization,
alongside the increase in testing facilities to identify more
individuals with COVID-19. Their investigation demonstrated that the
integration of these two control measures proved highly effective in
curbing the spread of the virus. The model by Srivastav et al.
\cite{Ref11} is given below:
\begin{equation}\label{Ref11}
\begin{aligned}
\frac{d S_1}{d t} & =-\beta_1 \left(I_1+I_2\right)S_1-\beta_3 S_1 H \\
\frac{d E_1}{d t} & =\beta_1 \left(I_1+I_2\right)S_1+\beta_3 S_1 H-\gamma_1 E_1\left(I_1+I_2\right)-\eta_1 E_1 \\
\frac{d I_1}{d t} & =\eta_1 E_1+\gamma_1 E_1\left(I_1+I_2\right)-\nu_1 I_1-\delta_1 I_1 \\
\frac{d S_2}{d t} & =-\beta_2 S_2\left(I_1+I_2\right)-\beta_4 S_2 H \\
\frac{d E_2}{d t} & =\beta_2 S_2\left(I_1+I_2\right)+\beta_4 S_2 H-\gamma_2 E_2\left(I_1+I_2\right)-\eta_2 E_2 \\
\frac{d I_2}{d t} & =\eta_2 E_2+\gamma_2 E_2\left(I_1+I_2\right)-\delta_2 I_2-\nu_2 I_2 . \\
\frac{d H}{d t} & =\nu_1 I_1+\nu_2 I_2-\delta_3 H-\alpha H . \\
\frac{d R}{d t} & =\alpha H .
\end{aligned}
\end{equation}

Where the model variables are defined as follows:  $S_1$, $E_1$,
$I_1$ are susceptible, exposed, and infected young individuals,
while
  $S_2$, $E_2$, $I_2$ represent susceptible, exposed, and infected old individuals. Home-isolated/hospitalized infected and recovered individuals of both groups are represented by $H$ and $R$ respectively.  The parameters are defined in Table \ref{Ref11Table} in the Appendix section.

\subsection{Analyzing Cost-Effective Control Measures for Co-Infections of Dengue and COVID-19 in Brazil}
Omame \textit{et al.} \cite{Ref22} explored a mathematical model for
COVID-19 and dengue, incorporating optimal control and
cost-effectiveness analysis. The local asymptotic stability of the
model was assessed under the condition that reproduction numbers
remain below unity. Utilizing data collected from February 1, 2021,
to September 20, 2021, the researchers fitted the model to
cumulative confirmed daily COVID-19 cases and deaths in Brazil. The
investigation included the estimation of several parameters,
including the rates of COVID-19 transmission and mortality, as well
as the attenuation of immunity acquired from infection. The study
extends the model to incorporate optimal control strategies,
establishing conditions for their existence and deriving the
optimality system using Pontryagin’s Principle. Simulation results
highlight the significant impact of dengue-only or COVID-19-only
control strategies in reducing new co-infection cases. Furthermore,
the findings emphasize the effectiveness of dengue prevention
strategies in averting new COVID-19 infections and highlight the
cost-effectiveness of controlling incident dengue infections to
mitigate co-infections.

 The model studied by Omame et al \cite{Ref22} is given below:
    \begin{equation}\label{Ref22}
    \begin{split}
    \frac{dS_{h}}{dt} & = \omega_{h} - \left(\frac{\Lambda_{vd}I_{vd}}{N_{h}} + \frac{\Lambda_{hc}(I_{hc} + I_{dc})}{N_{h}} \right) S_{h}- \varrho_{h} S_{h} + \eta_{hd}R_{hd} + \eta_{hc}R_{hc}, \\
    \frac{dI_{hd}}{dt} & = \frac{\Lambda_{vd}I_{vd}}{N_{h}} (S_{h} + R_{hc}) -(\alpha_{hd} +\varrho_{h}  + \varphi_{hd}) I_{hd} - \vartheta_{1}\frac{\Lambda_{hc}(I_{hc} + I_{dc})}{N_{h}}I_{hd} + \alpha_{hc}I_{dc},\\
    \frac{dR_{hd}}{dt} & = \alpha_{hd} I_{hd} - \varrho_{h} R_{hd} - \eta_{hd} R_{hd} - \frac{\Lambda_{hc}(I_{hc} + I_{dc})}{N_{h}} R_{hd},\\
    \frac{dI_{hc}}{dt} & = \frac{\Lambda_{hc}(I_{hc} + I_{dc})}{N_{h}} (S_{h} + R_{hd}) - (\alpha_{hc} + \varrho_{h}  + \varphi_{hc}) I_{hc} - \vartheta_{2}\frac{\Lambda_{vd}I_{vd}}{N_{h}}I_{hc} + \alpha_{hd}I_{dc}, \\
    \frac{dR_{hc}}{dt} & = \alpha_{hc} I_{hc} - \varrho_{h} R_{hc} - \eta_{hc} R_{hc} - \frac{\Lambda_{vd}I_{vd}}{N_{h}} R_{hc},\\
    \frac{dI_{dc}}{dt} & = \vartheta_{1}\frac{\Lambda_{hc}(I_{hc} + I_{dc})}{N_{h}}I_{hd} + \vartheta_{2}\frac{\Lambda_{vd}I_{vd}}{N_{h}}I_{hc} - (\varrho_{h} + \varphi_{hd} + \varphi_{hc} + \alpha_{hd} + \alpha_{hc})I_{dc}, \\
    \frac{dS_{vd}}{dt} & = \omega_{d} - \frac{\Lambda_{hd}(I_{hd} + I_{dc})}{N_{h}}S_{vd} - \varrho_{v}S_{vd}, \\
    \frac{dI_{vd}}{dt} & = \frac{\Lambda_{hd}(I_{hd} + I_{dc})}{N_{h}}S_{vd} - \varrho_{v} I_{vd}, \\
    \end{split}
    \end{equation}
The model variables are defined as follows: $S_{h}$ for susceptible
humans, $I_{hd}$ and $R_{hd}$ for infectious and recovered humans
from dengue, $I_{hc}$ and $R_{hc}$ for infectious and recovered
humans from COVID-19. Additionally, $I_{dc}$, $S_{v}$, and $I_{vd}$
represent infectious individuals co-infected with dengue and
COVID-19, susceptible vectors, and infectious vectors with dengue,
respectively. The parameters are defined in Table \ref{Ref22Table}
in the Appendix section.

\subsection{Analyzing the co-dynamics of COVID-19 and Malaria}
Tchoumi \textit{et al.} \cite{Ref23} developed a mathematical model
to investigate the co-dynamics of malaria and COVID-19, particularly
addressing the challenges in tropical and subtropical regions. The
stability conditions for equilibria in malaria-only and
COVID-19-only sub-models were established, indicating global
asymptotic stability for COVID-19 and potential backward bifurcation
in malaria under specific conditions. The study explored teh optimal
control strategies using Pontryagin's Maximum Principle to mitigate
the spread of both diseases. Simulation results indicated that
concurrently applying preventive measures for malaria and COVID-19
is more effective than implementing individual measures. While not a
case study, the research emphasized the importance of future
investigations considering therapeutic strategies, immunity
acquisition, treatment efficacy, and cost-effectiveness in disease
control. The mathematical model investigated by Tchoumi \textit{et
al.} \cite{Ref23}  is given as follows:

\begin{equation}\label{Ref23}
\left\{\begin{array}{l}
\dot{S}_h=\Lambda_h+\omega_m I_m+\omega_c I_c+\omega_{m c} I_{m c}-\left(\lambda_m+\lambda_c+\mu\right) S_h, \\
\dot{E}_m=\lambda_m S_h-\left(\lambda_c+\phi_m+\mu\right) E_m, \\
\dot{E}_c=\lambda_c S_h-\left(\lambda_m+\phi_c+\mu\right) E_c, \\
\dot{E}_{m c}=\lambda_c E_m+\lambda_m E_c-\left(\phi_{m c}+\mu\right) E_{m c}, \\
\dot{I}_m=\phi_m E_m-\left(\delta \lambda_c+\omega_m+\mu\right) I_m, \\
\dot{I}_c=\phi_c E_c-\left(\epsilon \lambda_m+\omega_c+\mu\right) I_c, \\
\dot{I}_{m E_c}=\delta \lambda_c I_m-(\sigma+\mu) I_{m E_c}, \\
\dot{I}_{m c}=\sigma I_{m E_c}+\phi_{m c} E_{m c}+\gamma I_{c E_m}-\left(\omega_{m c}+\mu\right) I_{m c}, \\
\dot{I}_{c E_m}=\epsilon \lambda_m I_c-(\gamma+\mu) I_{c E_m}, \\
\dot{S}_v=\Lambda_v-\left(\lambda_v+\mu_v\right) S_v, \\
\dot{E}_v=\lambda_v S_v-\left(\phi_v+\mu_v\right) E_v, \\
\dot{I}_v=\phi_v E_v-\mu_v I_v
\end{array}\right.
\end{equation}
The model variables are defined as follows: $S_{h}$ represents
susceptible individuals, $E_{m}$ and $I_{m}$ denote individuals who
are exposed to and infected with malaria only, respectively.
Similarly, $E_{c}$ and $I_{c}$ refer to individuals who are exposed
to and infectious with COVID-19. Additionally, the variables $E_{m
c}$ and $I_{m c}$ represent individuals who are exposed to and
infected with both malaria and COVID-19. Furthermore, the variables
$I_{m E_c}$ and $I_{c E_m}$ represent individuals infected with
malaria and exposed to COVID-19, and individuals infected with
COVID-19 and exposed to malaria, respectively. The parameters are
defined in Table \ref{TchoumiTable} in the Appendix section.
\subsection{Optimal Control Measures: Co-circulation of COVID-19 and Arboviruses in Espirito Santo, Brazil}
Omame et al.  \cite{Ref24} investigated the co-circulation of
COVID-19 and arboviruses (dengue, chikungunya, and Zika) using a
mathematical model. To calibrate their mathematical model, authors
utilized recorded data on COVID-19 and arboviruses from Espirito
Santo, a city in Brazil. The study employed optimal control
measures, including time-dependent interventions, to minimize
infectious cases and associated costs. The findings demonstrated
that targeted COVID-19 prevention measures in Espirito Santo
significantly reduced co-infections with Zika, Dengue, and
Chikungunya.However, the study emphasized that a sole focus on
COVID-19 prevention might have limited impact on arboviruses,
highlighting the necessity of integrated control strategies.   The
research provides valuable insights into disease interactions and
emphasizes the importance of comprehensive preventive strategies in
regions with the co-circulation of these infectious diseases.
Simulations suggested that a combined strategy against COVID-19,
Zika, Dengue, and Chikungunya is the most effective in minimizing
their co-circulation within a community.

\clearpage
\newpage
The model in \cite{Ref24} is given below:
\begin{equation}\label{Ref24}
\begin{split}
\frac{d\mathcal{S}^h}{dt} & = \Psi^h  -\left(\frac{\beta_{1} \mathcal{I}_{C}^h}{\mathcal{N}^h} + \frac{\beta_{2}(\mathcal{I}_{Z}^h + \mathcal{I}_{CZ}^h)+\beta^h_{2}\mathcal{I}_{Z}^v}{\mathcal{N}^h} + \frac{\beta^h_{3} \mathcal{I}_{D}^v}{\mathcal{N}^h} + \frac{\beta^h_{4} \mathcal{I}_{K}^v}{\mathcal{N}^h} + \vartheta^h  \right)\mathcal{S}^h\\
\frac{d\mathcal{I}_{C}^h}{dt}  & = \frac{\beta_{1}\mathcal{I}_{C}^h}{\mathcal{N}^h} (\mathcal{S}^h + \mathcal{R}_{C}^h + \mathcal{R}_{Z}^h + \mathcal{R}_{D}^h + \mathcal{R}_{K}^h)  -\left(\eta_{C} + \zeta_{C} + \vartheta^h  \right)\mathcal{I}_{C}^h \\
& - \frac{\beta_{2}(\mathcal{I}_{Z}^h + \mathcal{I}_{CZ}^h)+\beta^h_{2}\mathcal{I}_{Z}^v}{\mathcal{N}^h} \mathcal{I}_{C}^h - \frac{\beta^h_{3}\mathcal{I}_{D}^v}{\mathcal{N}^h} \mathcal{I}_{C}^h - \frac{\beta^h_{4}\mathcal{I}_{K}^v}{\mathcal{N}^h} \mathcal{I}_{C}^h + \zeta_{Z} \mathcal{I}_{CZ}^h + \zeta_{D} \mathcal{I}_{CD}^h + \zeta_{K} \mathcal{I}_{CK}^h\\
\frac{d\mathcal{I}_{Z}^h}{dt}  & = \frac{\beta_{2}(\mathcal{I}_{Z}^h + \mathcal{I}_{CZ}^h)+\beta^h_{2}\mathcal{I}_{Z}^v}{\mathcal{N}^h} (\mathcal{S}^h + \mathcal{R}_{C}^h + \mathcal{R}_{Z}^h + \mathcal{R}_{D}^h + \mathcal{R}_{K}^h) -\left(\eta_{Z} + \zeta _{Z} +\vartheta^h  \right)\mathcal{I}_{Z}^h \\
& -\frac{\beta_{1}\mathcal{I}_{C}^h}{\mathcal{N}^h} \mathcal{I}_{Z}^h + \zeta_{C} \mathcal{I}_{CZ}^h\\
\frac{d\mathcal{I}_{D}^h}{dt}  & = \frac{\beta^h_{3}\mathcal{I}_{D}^v}{\mathcal{N}^h} (\mathcal{S}^h + \mathcal{R}_{C}^h + \mathcal{R}_{Z}^h + \mathcal{R}_{D}^h + \mathcal{R}_{K}^h) -\left(\eta_{D} + \zeta _{D} +\vartheta^h  \right)\mathcal{I}_{D}^h - \frac{\beta_{1}\mathcal{I}_{C}^h}{\mathcal{N}^h} \mathcal{I}_{D}^h + \zeta_{C} \mathcal{I}_{CD}^h\\
\frac{d\mathcal{I}_{K}^h}{dt}  & = \frac{\beta^h_{4}\mathcal{I}_{K}^v}{\mathcal{N}^h} (\mathcal{S}^h + \mathcal{R}_{C}^h + \mathcal{R}_{Z}^h + \mathcal{R}_{D}^h + \mathcal{R}_{K}^h) -\left(\eta_{K} + \zeta _{K} +\vartheta^h  \right)\mathcal{I}_{K}^h - \frac{\beta_{1}\mathcal{I}_{C}^h}{\mathcal{N}^h} \mathcal{I}_{K}^h + \zeta_{C} \mathcal{I}_{CK}^h\\
\frac{d\mathcal{I}_{CZ}^h}{dt}  & = \frac{\beta_{2}(\mathcal{I}_{Z}^h + \mathcal{I}_{CZ}^h)+\beta^h_{2}\mathcal{I}_{Z}^v}{\mathcal{N}^h} \mathcal{I}_{C}^h + \frac{\beta_{1}\mathcal{I}_{C}^h}{\mathcal{N}^h} \mathcal{I}_{Z}^h -\left(\eta_{C} + \eta_{Z} + \zeta _{C} + \zeta _{Z} +\vartheta^h  \right)\mathcal{I}_{CZ}^h \\
\frac{d\mathcal{I}_{CD}^h}{dt}  & = \frac{\beta^h_{3}\mathcal{I}_{D}^v}{\mathcal{N}^h} \mathcal{I}_{C}^h + \frac{\beta_{1}\mathcal{I}_{C}^h}{\mathcal{N}^h} \mathcal{I}_{D}^h -\left(\eta_{C} + \eta_{D} + \zeta _{C} + \zeta _{D} +\vartheta^h  \right)\mathcal{I}_{CD}^h \\
\frac{d\mathcal{I}_{CK}^h}{dt}  & = \frac{\beta^h_{4}\mathcal{I}_{K}^v}{\mathcal{N}^h} \mathcal{I}_{C}^h + \frac{\beta_{1}\mathcal{I}_{C}^h}{\mathcal{N}^h} \mathcal{I}_{K}^h -\left(\eta_{C} + \eta_{K} + \zeta _{C} + \zeta _{K} +\vartheta^h  \right)\mathcal{I}_{CK}^h \\
\frac{d\mathcal{R}_{C}^h}{dt}  & = \zeta_{C} \mathcal{I}_{C}^h  -\left(\vartheta^h + \frac{\beta_{1}\mathcal{I}_{C}^h}{\mathcal{N}^h} + \frac{\beta_{2}(\mathcal{I}_{Z}^h + \mathcal{I}_{CZ}^h)+\beta^h_{2}\mathcal{I}_{Z}^v}{\mathcal{N}^h} +\frac{\beta^h_{3} \mathcal{I}_{D}^v}{\mathcal{N}^h} + \frac{\beta^h_{4} \mathcal{I}_{K}^v}{\mathcal{N}^h}\right)\mathcal{R}_{C}^h\\
\frac{d\mathcal{R}_{Z}^h}{dt}  & = \zeta_{Z} \mathcal{I}_{Z}^h -\left(\vartheta^h + \frac{\beta_{1}\mathcal{I}_{C}^h}{\mathcal{N}^h} + \frac{\beta_{2}(\mathcal{I}_{Z}^h + \mathcal{I}_{CZ}^h)+\beta^h_{2}\mathcal{I}_{Z}^v}{\mathcal{N}^h} +\frac{\beta^h_{3} \mathcal{I}_{D}^v}{\mathcal{N}^h} + \frac{\beta^h_{4} \mathcal{I}_{K}^v}{\mathcal{N}^h} \right)\mathcal{R}_{Z}^h\\
\frac{d\mathcal{R}_{D}^h}{dt}  & = \zeta_{D} \mathcal{I}_{D}^h -\left(\vartheta^h + \frac{\beta_{1}\mathcal{I}_{C}^h}{\mathcal{N}^h} + \frac{\beta_{2}(\mathcal{I}_{Z}^h + \mathcal{I}_{CZ}^h)+\beta^h_{2}\mathcal{I}_{Z}^v}{\mathcal{N}^h} +\frac{\beta^h_{3} \mathcal{I}_{D}^v}{\mathcal{N}^h} + \frac{\beta^h_{4} \mathcal{I}_{K}^v}{\mathcal{N}^h} \right)\mathcal{R}_{D}^h\\
\frac{d\mathcal{R}_{K}^h}{dt}  & = \zeta_{K} \mathcal{I}_{K}^h -\left(\vartheta^h + \frac{\beta_{1}\mathcal{I}_{C}^h}{\mathcal{N}^h} + \frac{\beta_{2}(\mathcal{I}_{Z}^h + \mathcal{I}_{CZ}^h)+\beta^h_{2}\mathcal{I}_{Z}^v}{\mathcal{N}^h} +\frac{\beta^h_{3} \mathcal{I}_{D}^v}{\mathcal{N}^h} + \frac{\beta^h_{4} \mathcal{I}_{K}^v}{\mathcal{N}^h}\right)\mathcal{R}_{K}^h\\
\frac{d\mathcal{S}^v}{dt} & = \Psi^v  -\left(\frac{\beta^v_{2}(\mathcal{I}_{Z}^h + \mathcal{I}_{CZ}^h)}{\mathcal{N}^h} + \frac{\beta^v_{3} (\mathcal{I}_{D}^h + \mathcal{I}_{CD}^h)}{\mathcal{N}^h} + \frac{\beta^v_{4} (\mathcal{I}_{K}^h + \mathcal{I}_{CK}^h)}{\mathcal{N}^h} + \vartheta^v  \right)\mathcal{S}^v \\
\frac{d\mathcal{I}_{Z}^v}{dt}  & = \frac{\beta^v_{2}(\mathcal{I}_{Z}^h + \mathcal{I}_{CZ}^h)}{\mathcal{N}^h} \mathcal{S}^v - \vartheta^v \mathcal{I}_{Z}^v \\
\frac{d\mathcal{I}_{D}^v}{dt}  & = \frac{\beta^v_{3}(\mathcal{I}_{D}^h + \mathcal{I}_{CD}^h)}{\mathcal{N}^h} \mathcal{S}^v - \vartheta^v \mathcal{I}_{D}^v \\
\frac{d\mathcal{I}_{K}^v}{dt}  & = \frac{\beta^v_{4}(\mathcal{I}_{K}^h + \mathcal{I}_{CK}^h)}{\mathcal{N}^h} \mathcal{S}^v - \vartheta^v \mathcal{I}_{K}^v \\
\end{split}
\end{equation}
Where the model variables are defined as follows:   $\mathcal{S}^h
$, $\mathcal{I}_{C}^h $, $\mathcal{I}_{Z}^h$,
$\mathcal{I}_{D}^h$,$\mathcal{I}_{K}^h$,   represent susceptible
humans, infectious humans with COVID-19, Zika virus, Dengue virus,
and Chikungunya, respectively. $\mathcal{I}_{CZ}^h $,
$\mathcal{I}_{CD}^h $,  $\mathcal{I}_{CK}^h $ represent humans
co-infected with Zika virus and COVID-19, co-infected with  Dengue
virus and COVID-19, and co-infected with  Chikungunya and COVID-19,
respectively.  $\mathcal{R}_{C}^h , \mathcal{I}_{Z}^h ,
\mathcal{R}_{D}^h, \mathcal{R}_{K}^h $ represent humans who have
recovered from COVID-19, Zika virus, Dengue virus, and Chikungunya,
respectively. For the mosquito population,   $\mathcal{S}^v $
represents susceptible mosquitoes, while $\mathcal{I}_{Z}^v $,
$\mathcal{I}_{D}^v $, $\mathcal{I}_{K}^v $ represent mosquitoes
infected with Zika virus, Dengue virus, and Chikungunya,
respectively. The model parameters are defined in Table
\ref{Ref24Table} in the Appendix section.

\subsection{A cost-effectiveness Analysis of COVID-19 Vaccination Strategies for Turkey}
Hagens et al. (2021) \cite{Ref14} conducted a study exploring the
economic implications of COVID-19 vaccination strategies in Turkey.
Their research, as of March 2021, focused on estimating the
cost-effectiveness of vaccination in comparison to a baseline
without vaccination and implemented measures. Employing an enhanced
SIRD model, the researchers considered different scenarios for the
initial year following vaccination. The findings revealed that, in
the context of healthcare, COVID-19 vaccination in Turkey is
cost-effective. The incremental cost-effectiveness ratios (ICER)
were observed to be 511 USD/QALY and 1045 USD/QALY under different
scenarios of vaccine effectiveness. The societal perspective
revealed cost savings in both scenarios, emphasizing that a minimum
vaccine uptake of at least 30\% is required for cost-effectiveness.
Sensitivity and scenario analyses, along with iso-ICER curves,
confirmed the robustness of the findings. The article proposed that
COVID-19 vaccination in Turkey is highly cost-effective, potentially
cost-saving, particularly when considering the economic consequences
of potential lockdowns in the absence of access to vaccination. The
following differential equations considered in the model by Hagens
\textit{et al.} \cite{Ref14}.
        \begin{equation}\label{Ref14}
        \begin{gathered}
        \frac{d S_i}{d t}=\left(-\beta S_i \sum_j C_{i j} I_j / N_j\right)-V_i \\
        \frac{d R_i}{d t}=f_i\left(F H m_i\left(H m_i\right)+F I c_i\left(I c_i\right)+F H s_i\left(H s_i\right)\right)+V_i
        \end{gathered}
        \end{equation}

        In the given context, $V_i$ represents the count of vaccinated individuals, while $F H m$, $Fic$, and $F H s$ denote the inverses of the durations for staying at home, being in the ICU, and hospitalization if unwell, respectively. Additionally, $\mathrm{Hm}$, $Ic$, and $Hs$ represent the count of individuals experiencing illness at home, in the ICU, and in hospitals, respectively and $f$ signifies the rate of recovery.

    \subsection{The cost-effectiveness and epidemiological impact of COVID-19 vaccination across diverse
regions} We delve into a series of research papers, each
contributing valuable insights into the cost-effectiveness and
epidemiological impact of COVID-19 vaccination across diverse
regions. We have provided the main findings of these studies without
going into much detail about the models and methods used in them.
    \subsubsection{A Case Study from Catalonia}
    Utilizing data from Catalonia, the authors \cite{Ref12} constructed and evaluated  a cost-effective COVID-19 immunization model. Findings suggest that not only is the mass vaccination campaign cost-saving, but it also proves to be widely cost-effective for the health system. The efficient allocation of resources leads to substantial social and economic advantages.
    \subsubsection{A Case Study from Kenya}
    The authors \cite{Ref13} conducted a 1.5-year evaluation of COVID-19 vaccine cost-effectiveness in Kenya from a societal standpoint.   Employing an age-structured epidemic model, they assumed a baseline of natural immunity in at least 80\% of the population before the emergence of the immunological escape variety.
   The study explored the effects of vaccine coverage (30\%, 50\%, or 70\%) among the adult population ($>18$) through both slow and rapid roll-out scenarios, prioritizing individuals over 50 (80\% uptake in all scenarios). Cost data, obtained from initial analyses, indicated vaccine delivery costs ranging from US\$3.90 to US\$6.11 per dosage and vaccine procurement costs of US\$7 per dose. The study concluded that immunizing young adults might not be cost-effective due to earlier exposure and limited protection across the majority of the Kenyan population.
\subsubsection{A Case Study from Korea}

The authors \cite{Ref15} conducted a comprehensive
cost-effectiveness analysis of oral antivirals targeting SARS-CoV-2
infection in Korea, utilizing an epidemic model. Their research
demonstrated that the use of oral treatment, specifically
nirmatrelvir/ritonavir, among older patients with symptomatic
COVID-19, emerged as the most cost-effective approach. Moreover,
this therapeutic strategy was associated with a substantial
reduction in the demand for new hospital admissions in Korea.

\subsubsection{A Case Study from the United States}

In the United States, Li et al. \cite{Ref16} investigated the
cost-effectiveness of BNT162b2 COVID-19 booster doses. Their
research concluded that administering BNT162b2 booster doses to
senior Americans (aged $\ge$ 65 years) represented the most
economical option. Furthermore, the authors observed that in regions
with high viral transmission rates, even less effective COVID-19
vaccinations and booster doses might still be deemed economically
viable.

\subsubsection{Two Case Studies from Hong Kong}

    Using Hong Kong as a case study, Tao et al \cite{Ref17} conducted research on COVID-19 outbreak predictions and cost-effectiveness analysis. Results from their analysis indicated that the total number of patients in Hong Kong with new infections may be decreased by 83.89\% by strictly enforcing quarantine measures in the third stage of COVID-19. In conjunction with the vaccine plan, vaccinating 90\% of the population within a set period of time could significantly reduce the epidemic in Hong Kong. Additionally, there may be financial gains and a 10.74 percent cost reduction compared to the non-vaccination situation.

Xiong et al. \cite{hong} studied the economic impact of COVID-19
vaccination programs in Hong Kong. They used a Markov model with a
susceptible–infected–recovered structure over a 1-year time horizon.
The research revealed that despite a high initial cost, the
vaccination program became economically justified during periods of
increased infection rates, like the Omicron wave. The study
emphasized the program's effectiveness in reducing infection and
mortality rates. However, it also acknowledged the substantial
economic burden, highlighting the need for strategic implementation
and prioritization. The study concluded that the focus should be on
the elderly population to improve vaccine coverage.

\subsubsection{A Case Study from Six U.S. Cities}
Zang et al. \cite{Ref18} explored the cost-effectiveness of linked,
opt-out HIV testing. Their study also explored the potential
epidemiological impact of COVID-19 on HIV epidemic, employing a case
study across six U.S. cities. The study's findings revealed that in
the absence of linked, opt-out HIV testing, a worst-case
scenario—characterized by no behavioral change and a 50\% reduction
in service access—could lead to an estimated 9.0 percent increase in
HIV infections. Conversely, in the best-case scenario—entailing a
50\% reduction in risk behaviors and no service disruptions—the
study anticipated a substantial estimated decrease of 16.5 percent
in HIV infections between 2020 and 2025. Additionally, the authors
predicted that implementing HIV testing at varying levels (ranging
from 10 to 90 percent) might prevent a total of 576 to 7,225 new
infections.

\subsubsection{A Case Study from South Africa}

    Reddy \textit{et al.} \cite{Ref19} investigated the impact of COVID-19 vaccination in South Africa. Their research findings demonstrate that ensuring vaccine coverage for at least 40\% of the population and prioritizing efficient vaccine rollout resulted in a substantial preventive effect. This approach successfully averted more than 9 million new COVID-19 infections, prevented over 70,000 deaths, and contributed to cost reductions by minimizing the need for hospitalizations.

\subsubsection{A Case Study from Low- and Middle-Income Countries}
Utami \textit{et al.} \cite{Ref20} investigated the
cost-effectiveness of COVID-19 vaccination in low- and middle-income
nations. Their research revealed that administering COVID-19
vaccines proved to be both cost-effective and cost-saving in
reducing mortality and limiting the spread of COVID-19 within low-
and middle-income countries.

    \subsubsection{A Case Study from China}

Zhou \textit{et al.} \cite{Ref21} examined the cost-effectiveness of
COVID-19 vaccination in China. According to their research,
delivering the COVID-19 vaccine to the general population in
mainland China might reduce the infection rate by 55\% and the death
rate by 3.7\% compared to a no-vaccination scenario.

    This comprehensive overview of diverse studies provides a nuanced understanding of the global landscape of COVID-19 vaccination strategies. Each case study contributes unique insights into the economic considerations and epidemiological outcomes of vaccination efforts across different regions.

    \section{Discussion and future directions}\label{Discussion-concluding remarks}
    In this paper, we conducted an extensive literature review focusing on articles related to the cost-effectiveness of COVID-19 vaccination. Notably, various authors have developed diverse epidemic models to assess the cost-effectiveness of COVID-19 vaccination. The majority of research findings consistently emphasize the effectiveness of COVID-19 vaccination and other preventive measures as the most cost-effective strategies against the disease. While the accumulated evidence supports the significance of enhanced vaccination and preventive strategies, it is crucial to acknowledge a significant hurdle in the form of vaccine nationalism, particularly affecting low and middle-income economies. This remains a major obstacle in realizing an effective and widespread vaccination program. Consequently, further research is warranted to delve into how vaccine nationalism might influence or hinder the success of a cost-effective vaccination program for COVID-19 and other vaccine-preventable diseases. The authors strongly advocate for a more comprehensive exploration of these dynamics to inform more inclusive and globally impactful vaccination strategies.

The economic and social impact of any pandemic cannot be
overemphasized. Infectious diseases could destabilize the growth and
progress of every economy. Reduced human capital is a consequential
outcome, stemming from both a decline in labor supply and diminished
productivity caused by sickness and mortality resulting from the
infectious disease. In this review, we have presented detailed
reports on COVID-19 epidemic models, specifically focusing on the
cost-effectiveness of COVID-19 vaccination. Additionally, we have
identified areas in the existing literature that require further
investigation and improvement. These areas represent potential gaps
or limitations in current research, and addressing them can
contribute to a more comprehensive understanding of the dynamics
surrounding pandemics and vaccination strategies. Moreover, we have
proposed directions for future research to fortify our preparedness
against the occurrence of future pandemics.

Despite the significant success in our research, notable gaps
persist in our knowledge and understanding of vaccine nationalism
and macroeconomic epidemiological models for COVID-19, warranting
further investigation. Specifically, future research should offer a
detailed review of COVID-19 vaccine nationalism, shedding light on
the uneven distribution of COVID-19 vaccines between wealthy and
impoverished nations during the initial rollout. Additionally, a
comprehensive review of macroeconomic epidemiological models is
essential to explore how these models can provide effective
strategies for managing the global economy in the event of future
pandemics.

Given the current state of research and the prevailing trends,
future directions in the study of the cost-effectiveness of COVID-19
vaccination should explore incorporating diffusive models that
consider both temporal and spatial dimensions. The existing body of
literature primarily emphasizes epidemic models highlighting the
efficacy of vaccination and preventive measures. However, there is a
need to elevate the complexity of these models by integrating
diffusion dynamics across both time and space. By infusing diffusive
models, researchers can attain a more detailed grasp of how
infectious diseases propagate and the ramifications of vaccination
strategies across diverse geographical locations and time frames.
This approach has the potential to furnish more precise predictions
and insights into the cost-effectiveness of vaccination programs,
particularly in light of the challenges presented by vaccine
nationalism. By incorporating spatial and temporal dimensions,
future research can contribute to the development of more adaptable
and robust strategies that address the complexities of global health
crises.
\newpage

\newpage
%\reftitle{References}
%\begin{adjustwidth}{-\extralength}{0cm}

%\end{adjustwidth}
\newpage

\section*{Appendix}
\begin{table}[h!]
    \begin{tabular}{ll}
        \hline
        \textbf{Parameter} & \textbf{Description} \\
        \hline
        $B(t)$ & Concentration of the SARS-CoV-2 in the environment. \\
        $\beta_1, \beta_2, \beta_3$ & Transmission rates by individuals \\
        & in the infected classes $E(t), A(t), I(t)$, respectively\\
        $\beta_4$ & Propensity rate of susceptible individuals\\
        &  getting the virus through the environment. \\
        $(1-\tau)\delta$ & Rate at which the exposed individuals\\
        &  develops symptoms\\
        $\tau \delta$ & Rate of new asymptomatic infection\\
        $d_1$ & Disease-induced death rate\\
        $\gamma_1$ and $\gamma_2$ & recovery rates for symptomatic and asymptomatic\\
        $\psi_1, \psi_2$ and $\psi_3$& Virus shedding rates into the environment by \\
        & the exposed, infected and asymptomatically infected\\
        $\phi$ & natural removal of the virus from the environment\\
        \hline
        \end{tabular}
        \caption{Description of the parameters of the model \eqref{Ref5}}
        \label{Ref5Table}
        \end{table}

\begin{table}[h!]
    \begin{tabular}{ll}
        \hline
        \textbf{Parameter} & \textbf{Description} \\
        \hline
        $\Lambda$ & Recruitment rate \\
        $\lambda$ & Force of infection \\
        $\beta$ & Natural death rate \\
        $\beta_1$ & Transmission rate \\
        $\eta$ & Transmission rate from the environment \\
        $\gamma$ & Relative transmissibility of asymptomatic individuals \\
        $k_1$ & Proportion of individuals who are timely diagnosis \\
        $k_2$ & Progression rate from exposed to the symptomatic \\
        & (severely infected) class \\
        $\alpha$ & Progression rate from exposed to the asymptomatic class \\
        $\phi$ & Disease induced death rate \\
        $v_1$ & Proportion of asymptomatic patients who later move to \\
        & the symptomatic (severely infected) class \\
        $v_2$ & Progression from asymptomatic to the symptomatic \\
        & (severely infected) class \\
        $\epsilon$ & Progression from asymptomatic to the recovered class \\
        $\rho$ & The rate at which symptomatic (severely infected) \\
        & patients recovers \\
        $\tau_1$ & The rate at which the recovered individuals \\
        & join the susceptible class \\
        $m_1$ & Natural decay rate of virus from the environment (Surfaces) \\
        $m_2$ & The rate of viral release into the environment\\
        &  by asymptomatic patients \\
        \hline
    \end{tabular}
    \caption{Description of the parameters of the model \eqref{Ref6}}
    \label{Ref6Table}
\end{table}

    \begin{table}[h!]
        \begin{tabular}{ll}
            \hline
            \textbf{Parameter} & \textbf{Description} \\
            \hline
            $\beta$ & Effective transmission coefficient \\
            $\varepsilon_1$ & Modification parameter for a reduced \\
            & transmission from asymptomatic humans \\
            $\varepsilon_2$ & Modification parameter for a more reduced\\
            &  transmission from hospitalized class \\
            $\alpha$ & Rate of disease progression from exposed class \\
            $l_1$ & Proportion of exposed with symptoms after the incubation period \\
            $1-l_1$ & Proportion of exposed without symptoms after the incubation period \\
            $h_1$ & Hospitalization rate for symptomatic class \\
            $h_2$ & Hospitalization rate for asymptomatic class after confirmation \\
            $r_1$ & Recovery rate for symptomatic class \\
            $r_2$ & Recovery rate for asymptomatic class \\
            $\gamma$ & Recovery rate for hospitalized class \\
            $\delta_1$ & Disease-induced mortality rate for symptomatic class \\
            $\delta_2$ & Disease-induced mortality rate for hospitalized class \\
            \hline
        \end{tabular}
        \caption{Description of the parameters of the model \eqref{Ref7}}
        \label{Ref7Table}
    \end{table}

    \begin{table}[h!]
        \begin{tabular}{ll}
            \hline
            \textbf{Parameter} & \textbf{Description} \\
            \hline
            $\beta_1$ & The rate of people who were infected by contact with \\
            & infected and asymptomatic \\
            $\beta_2$ & The rate of people who were infected by contact with\\
            & infected and symptomatic \\
            $\alpha_1$ & The rate of exposed become infected and asymptomatic \\
            $\alpha_2$ & The rate of exposed become infected and symptomatic \\
            $\theta_1$ & The number of infected and asymptomatic who have been under lockdown \\
            $\theta_2$ & The number of infected and symptomatic who have been under lockdown \\
            $\delta_1$ & Mortality rate due to serious complications \\
            $\delta_2$ & mortality rate due to serious complications and with chronic diseases \\
            $\delta_3$ & The rate of people who died under quarantine in hospitals. \\
            $\mu$ & Natural mortality \\
            $\eta_1$ & The rate of infection with serious complications without chronic \\
            &diseases under lockdown \\
            $\eta_2$ & The rate of infection with serious complications with \\
            &chronic diseases under lockdown \\
            $\sigma$ & The rate at which people recover from the virus \\
            $\gamma$ & The rate at which infected become serious complications \\
            &and those without and with chronic disease \\
            $\chi$ & The fraction of infected and asymptomatic \\
            \hline
        \end{tabular}
        \caption{Description of the parameters of the model \eqref{Ref8}}
        \label{Ref8Table}
    \end{table}

\begin{table}[h!]
    \begin{tabular}{ll}
        \hline
        \textbf{Parameter} & \textbf{Description} \\
        \hline
        $\Lambda$ & Recruitment rate\\
        $\eta_1$ & Contact rate\\
        $\psi$ & transmissibility factor associated with the\\
        &  asymptomatically infected persons.\\
        $\mu$ & natural death rate\\
        $\theta$ & infection following exposure\\
        $\omega$ and $\rho$ & incubation periods of the disease\\
        $\tau_1, \tau_2, \phi_1, \phi_2$ & Recovery rates for infected,\\
        &  asymptomatically infected, quarantined, and hospitalized individuals\\
        $\gamma$ and $\delta_2$ & hospitalization rates for infected and quarantined individuals\\
        $\xi_1$ and $\xi_2$ & COVID-19-induced death rates\\
        $\eta_2$ & Contact rate due to a visit to the seafood market\\
        $q_1$ and $q_2$ & Coronaviruses shed rates\\
        $q_3$ & Rate of removal of coronaviruses from the seafood market\\
        $\delta_1$ & Quarantine rate\\
        \hline
    \end{tabular}
    \caption{Description of the parameters of the model \eqref{Ref9}}
    \label{Ref9Table}
\end{table}

    \begin{table}[h!]
        \centering
        \begin{tabular}{ll}
            \hline
            Parameter & Description \\
            \hline
            $\beta_1$ & Transmission rate from $I_1$ or $I_2$ to $S_1$, \\
            $\beta_2$ & Transmission rate from $I_1$ or $I_2$ to $S_2$, \\
            $\beta_3$ & Transmission rate from $H$ to $S_1$, \\
            $\beta_4$ & Transmission rate from $H$ to $S_2$ \\
            $\delta_1$ & Disease related death rate in $I_1$ compartment, \\
            $\delta_2$ & Disease related death rate in $I_2$ compartment, \\
            $\delta_3$ & Disease related death rate in $H$ compartment, \\
            $\nu_1$ & Rate of detection/isolation in $I_1$ compartment, \\
            $\nu_2$ & Rate of detection/isolation in $I_2$ compartment, \\
            $\eta_1$ & Rate of progression of individuals from $E_1$ to $I_1$, \\
            $\eta_2$ & Rate of progression of individuals from $E_1$ to $I_1$, \\
            $\gamma_1$ & Rate of reinfection in $E_1$ compartment, \\
            $\gamma_2$ & Rate of reinfection in $E_2$ compartment, \\
            $\alpha$ & Recovery rate of home isolated/hospitalized people. \\
            \hline
        \end{tabular}
        \caption{Description of the parameters of the model \eqref{Ref11}}
        \label{Ref11Table}
    \end{table}

    \begin{table}[h!]
        \centering
        \begin{tabular}{l l}
            \hline
            \textbf{Parameter} & \textbf{Interpretation}\\
            \hline
            $\omega_{h}$ & Human recruitment rate \\
            $\omega_{d}$ & Vector recruitment rate \\
            $\varrho_{h}$ & Human natural death rate \\
            $\eta_{hd}$ & Loss of infection acquired immunity to dengue \\
            $\Lambda_{vd}$ & Effective contact rate for vector to human transmission of dengue\\
            $\Lambda_{hd}$ & Effective contact rate for human to vector transmission of dengue\\
            $\alpha_{hd}$ & Dengue recovery rate\\
            $\Lambda_{hc}$ & Effective contact rate for human to human transmission of COVID-19\\
            $\eta_{hc}$ & Loss of infection acquired immunity to COVID-19\\
            $\varphi_{hc}$ & COVID-19-induced death rate\\
            $\alpha_{hc}$ & COVID-19 recovery rate \\
            $\vartheta_{1}$ & Modification parameter accounting for susceptibility of dengue-infected\\
            & Individuals to COVID-19\\
            $\vartheta_{2}$ & Modification parameter accounting for susceptibility of COVID-19-infected\\
            & Individuals to dengue\\
            $\varphi_{hd}$ & Dengue-induced death rate\\
            $\varrho_{v}$ & Vector removal rate\\
            \hline
        \end{tabular}
        \caption{Description of the parameters of the model \eqref{Ref22}}
        \label{Ref22Table}
    \end{table}

    \begin{table}[h!]
        \centering
        \begin{tabular}{l l}
            \hline
            \textbf{Parameter} & \textbf{Interpretation}\\
            \hline
            $\Lambda_h$ & Recruitment rate of the host population \\
            $\phi_m$ & Progression rate from exposed to infectious malaria state \\
            $\phi_c$ & Progression rate from exposed to infectious COVID-19 state \\
            $\phi_{m c}$ & Fraction of individuals moving to the co-infection $I_{m c}$ class \\
            $\omega_m$ & Recovery rate of malaria infected individuals \\
            $\omega_c$ & Recovery rate of COVID-19 infected individuals \\
            $\omega_{m c}$ & Recovery rate of malaria and COVID-19 infectious individuals \\
            $\mu$ & Death rate of the host population \\
            $\delta$ & Enhancement factor of acquiring COVID-19 following malaria infection \\
            $\epsilon$ & Enhancement factor of acquiring malaria following COVID-19 infection \\
            $\sigma$ & COVID-19 infection rate of individuals already infected with malaria \\
            $\gamma$ & Malaria infection rate of individuals already infected with COVID-19 \\
            $\kappa$ & Fraction of individuals employing personal protection \\
            $\zeta$ & Efficacy of personal protection \\
            $b$ & Number of female mosquito bites per day \\
            $\beta_m$ & Malaria transmission probability per mosquito bite \\
            $\beta_c$ & COVID-19 transmission probability per contact \\
            Mosquitoes & \\
            $\Lambda_v$ & Recruitment rate of vectors \\
            $\phi_v$ & Progression rate from exposed to infectious class \\
            $\mu_v$ & Natural death rate of vectors \\
            $\beta_v$ & Transmission probability in vectors from infected humans \\
            \hline
        \end{tabular}
        \caption{Description of the parameters of the model \eqref{Ref23} }
           \label{TchoumiTable}
    \end{table}

\begin{table}[h!]
    \begin{center}
        \begin{tabular}{l l}
            \hline
            \textbf{Parameter} & \textbf{Description}\\
            \hline
            $\Psi^v$  & vector recruitment rate \\
            $\Psi^h$  & human recruitment rate\\
            $\beta_1 $  & contact rate for COVID-19 infection\\
            $\beta_2$  &  zika infection contact rate (sexual transmission)\\
            $\beta^h_{2} $  & zika infection contact rate (vector human)\\
            $\beta^h_{3} $  & dengue infection contact rate (vector human)\\
            $\beta^h_{4} $  & chikungunya infection contact rate (vector human)\\
            $\beta^v_{2}$  & zika infection contact rate (human to vector)\\
            $ \beta^v_{3}$  & dengue infection contact rate (human to vector)\\
            $\beta^v_{4}$  & chikungunya infection contact rate (human to vector)\\
            $\vartheta^h $ & human natural death rate\\
            $\vartheta^v$ & vector removal rate\\
            $\eta_{C}, \eta_{Z}, \eta_{D}, \eta_{K}$ & COVID-19, zika, dengue and\\
            &  chikungunya disease-induced death rates, respectively\\
            $\zeta_{C}$  & COVID-19 recovery rate\\
            $\zeta_{Z}, \zeta_{D}, \zeta_{K}$  & zika, dengue, chikungunya recovery rates\\
            \hline
        \end{tabular}
    \end{center}
    \caption{Description of parameters of the model \eqref{Ref24}}
    \label{Ref24Table}
\end{table}

\end{document}